# FULL SPECTRAL TESTING OF LINEAR CONGRUENT METHOD WITH A MAXIMUM PERIOD

N. Temirgaliyev

**Abstract.** In this paper the spectral analysis of all possible linear congruent sequences with a maximum period is conducted and the best random number generators are selected among them.

## §1. Introduction

*A linear congruent sequence with a maximum period $N$ is referred to a recurrent sequence*

$$X_{n+1} = (aX_n + c) \bmod N \quad (n = 0,1,\ldots,N-1), \qquad (1.1)$$

which depends on the four «*magic*» positive integers $X_0, N, a, c$, where

$$
\begin{aligned}
&N, \text{ mod}ule && 0 < N \\
&a, \text{ multiplier} && 2 \leq a < N \\
&c, \text{ increment} && 1 \leq c < N, c \text{ and } N \text{ relatively prime} \\
&X_0, \text{ initial value} && 0 \leq X_0 < N
\end{aligned}
\qquad (1.2)
$$

such that

$$(a-1)^{\tau(a,N)} \equiv 0 \pmod{N} \quad \text{and} \quad (a-1)^{\tau(a,N)-1} \not\equiv 0 \pmod{N}, \qquad (1.3)$$

for $\tau(a,N) \geq 2$.

The indicator $\tau(a,N)$, uniquely determined by $a$ and $N$, is called a *potential* of the sequence (1.1)-(1.2). Further, for given $a$, $N$ and depending on them $\tau(a,N)$, a comparison (1.3) also uniquely determines an integer $\lambda(a,N) \equiv \dfrac{(a-1)^{\tau(a,N)}}{N}$ such that

$$(a-1)^{\tau(a,N)} = N\lambda(a,N) \quad (\tau(a,N) \geq 2, \quad 1 \leq \lambda(a,N) < (a-1)^{\tau(a,N)-1}). \qquad (1.4)$$

It can be emphasized and also required for future, that two potentials $\tau(a,N) \geq 2$ and $\lambda(a,N) \geq 1$ are uniquely connected with every $a$ and $N$.

By the definition of recurrence sequences (1.1) - (1.2), the first repetition of previously occurred number forms a cycle, which can be then repeated an infinite number of times. This repeating cycle is called *a period*.

N.Temirgaliyev, Professor, Director of the Institute of Theoretical Mathematics and Scientific Computation, L.N.Gumilyov Eurasian National University, 2 Satpayev Str., Astana, 010008 Republic of Kazakhstan, e-mail: ntmath10@mail.ru

It is well known (see [2, Section 3.2.1.3]), that the conditions (1.2)- (1.3) provide a necessary and sufficient condition, that a $N$-termed sequence (1.1) has a maximum period of length $N$.

This implies that in (1.1) each number $0,1,...,N-1$ appears only once, because in this case $X_0$ does not affect to the length of period, you can take any of these numbers for $X_0$. In this paper we consider that $X_0 = 0$.

Also $c, 1 \leq c < N$ is any positive integer number, which is mutually prime with $N$. It can be considered to be 1, that does not preclude the possibility of application to replace these two numbers on other permitted ones.

Following these agreements, the sequence (1.1) depends only on $a$ and $N$, $2 \leq a < N$, that can be reflected in the following notations

$$\langle X_n \rangle = \langle X_n(a,N,c,X_0) \rangle_{n=0}^{N-1} = \langle X_n(a,N) \rangle_{n=0}^{N-1} = \langle X_n(a,N) \rangle. \qquad (1.5)$$

From these notations (1.5), we often use the last (and sometimes the first short notation when it does not lead to confusions).

The last sequence(periodic with a period of $N$) (1.1)-(1.4) can be also called as a *random number generator*, or, more precisely, *Lehmer's random number generator* or sometimes in short just a *generator* or a *sequence* $\langle X_n \rangle$.

It is obvious thatthere are different requirements for randomness of the sequences (1.1)-(1.4), as well as for all other generators.

Numerous theoretical and empirical tests of randomness, at least the main ones are described in detail in [1-2].

All our attention is focused on the following ones [2, Part 3.3.4. A. **Ideas serving justification of criterion**]: *"The most important teststo check how random sequence is associated with the properties of the joint distributions of$^s$ successive elements of the sequence, and the spectral criterion is just used to test hypotheses about these distributions. If a sequence $U_n = \frac{X_n}{N}$ with a period $N$, then for the construction of criteria it is need to analyze the set of all $N$ points $\{(U_n, U_{n+1},...,U_{n+s-1}) | 0 \leq n < N\}$ in $s$-dimensional space."*

Thus, another parameter $s \geq 2$ is appeared, which is responsible for the independence of the sequence $s$-dimensional vectors

$$(X_n, X_{n+1},...,X_{n+(s-1)}),$$

in which the quantitative characteristic of independence is expressed in terms of the value of $s$-**dimensional accuracy of the random numbers generator** (the definition of $v_s$ is given in (2.5) - (2.7) in §2 below),

$$v_s(a,N) \equiv v_s(a,N; \text{ a connection between } a,N,s,\tau,\lambda \text{ in form of } (a-1)^\tau = N\lambda$$
$$\text{from (1.3)-(1.4))}, \qquad (1.6)$$

Thus, the task is to find the numbers $a, N$ and $s$ with conditions (1.1)-(1.4) and with a greater quantity value $v_s(a,N)$. The following theorem gives an upper value (see [1, 2, Section 3.3.4, E.].)

**Theorem B.** For all $a$, $N$ and $s$ the following inequality is true

$$v_s(a,N) \leq \gamma_s N^{\frac{1}{s}}, \qquad (1.7)$$

where $\gamma_s$ takes values

$$(4/3)^{1/4}, 2^{1/6}, 2^{1/4}, 2^{3/10}, (64/3)^{1/12}, 2^{3/7}, 2^{1/2}$$

for $s = 2,3,4,5,6,7,8$.

As any upper estimate, the inequality (1.7) can be greatly inflated, so the problem might be not solved (see [1, page 111.].) :*«As no one knows what are the best attainable $v_s$ value, it is difficult to determine exactly which $v_s$ values can be considered to be satisfied».*

Here, in Theorems 1-7 (see. Below §3), *«the best achievable values»* of $v_s$ *are found and given.*

It is worth pointing out that all possible sequences of (1.1) are not studied here, only those under the conditions of (1.2)-(1.4) with a maximum period $N$, which does not limit the generality, and cuts off less valuable.

Conditions (1.1)-(1.4) determine the relationship between $a$ and $N$ through $\tau(a,N)$ and $\lambda(a,N)$. Therefore, investigated problem with respect to the $s$ can be divided into the following disjoint in pairs cases

$$\begin{aligned} &1^0. \ s = \tau(a,N) = 2, \quad 1 \leq \lambda(a,N), \\ &2^0. \ s = \tau(a,N) \geq 3, \quad 1 \leq \lambda(a,N), \\ &3^0. \ 2 \leq s < \tau(a,N), \quad 1 \leq \lambda(a,N), \\ &4^0. \ 2 \leq \tau(a,N) < s, \quad 1 \leq \lambda(a,N). \end{aligned} \qquad (1.8)$$

Finally, we define the symbols. As always, $Z^s$ means the integer lattice in Euclidean space $R^s$, $Z_0^s$ is $Z^s$ without point $(0,...,0)$. A proposal *« A divide B »* is denoted by $A \mid B$.

For positive sequences $A_N$ and $B_N$ we introduce the notation $A_N \underset{\sim}{<} B_N$, if $A_N \leq \gamma_N B_N$, or what is the same $\gamma_N^{-1} A_N \leq B_N$, where $\gamma_N \to 1$, consequently $\gamma_N^{-1} \to 1$, for $N \to \infty$. By $A_N \approx B_N$ we determine a simultaneous execution of $A_N \underset{\sim}{<} B_N$ and $B_N \underset{\sim}{<} A_N$.

Throughout the paper we assume that $(t \geq 1) 1 < p_1 < \cdots < p_t$ ordinary numbers, - $\aleph_1,...,\aleph_t$ and $r_1,...,r_t$ positive integers, $N = p_1^{\aleph_1} \cdots p_t^{\aleph_t}$ and $a = d \times p_1^{r_1} \cdots p_t^{r_t} + 1$, where a positive integer $d$ is mutually prime with $N$, i.e., in the decomposition of $d$ into ordinary numbers $p_1,...,p_t$.

It can be noted that further in the formulation of the results the exact values $\gamma_N$ are specified, which provides effective computing applications of the results.

In the shortest summary, for explicitly given positive integers $\tau \geq 2, a \geq b_\tau + 1$ and $1 \leq \lambda < (a-1)^{\tau-1}$, the random number generator is the following

$$X_{n+1} = (aX_n + c) \bmod \frac{(a-1)^\tau}{\lambda} \left( n = 0,..., \frac{(a-1)^\tau}{\lambda} - 1 \right),$$

where $X_0 = 0$, and $1 \leq c < (a-1)^\tau / \lambda$ is any number, mutually prime with $a-1$, and this is the new results in spectral testing (ST), the following relations are conducted:

**ST-2:** $v_2^2\left(a, N; (a-1)^2 = N\right) = (a-1)^2 \left(1 - 2\frac{a-2}{(a-1)^2}\right) = N\left(1 - 2\frac{\sqrt{N}-1}{N}\right)$,

**ST** $(2 \leq s = \tau)$: $N^{\frac{2}{s}}\left(1 - (b_s - 1)N^{-\frac{1}{s}}\right)^2 = (a - b_s)^2 \leq v_s^2\left(a, N; (a-1)^s = N\right) \leq$

$$\leq a^2 + 1 = N^{\frac{2}{s}}\left(1 + 2N^{-\frac{1}{s}} + 2N^{-\frac{2}{s}}\right).$$

**ST** $(2 \leq s < \tau)$: $(N\lambda)^{\frac{2}{\tau}}\left(1 - (b_s - 1)(N\lambda)^{-\frac{1}{\tau}}\right)^2 = (a - b_\tau)^2 \leq v_s^2\left(a, N; (a-1)^\tau = N\lambda, 1 \leq \lambda \leq (a-1)^{\tau-s}\right) \leq$

$$\leq a^2 + 1 = (N\lambda)^{\frac{2}{\tau}}\left(1 + 2(N\lambda)^{-\frac{1}{\tau}} + 2(N\lambda)^{-\frac{2}{\tau}}\right),$$

**ST** $(s > \tau \geq 2)$: $v_s^2\left(a, N; (a-1)^\tau = N\lambda, \lambda \geq 1\right) \leq \sum_{k=0}^{\tau}\left(\binom{\tau}{k}\right)^2$.

By this, the problems discussed in *" C. The conclusion of computational method ([1, Section 3.3.4]): These examples illustrate the methods of use of the spectral test. However, there is a significant gap in our discussions: is there ever any possibility to determine the value v, not spending too much computer time? As an example, you can find out what exactly values $s_1 = 227$, $s_2 = 983$ и $s_3 = 130$ corresponds to the minimum amount of the sum $s_1^2 + s_2^2 + s_3^2$, following to the condition $s_1 + 3141592621 s_2 + 3141592621^2 s_3 \equiv 0 \pmod{10^{10}}$? It is obvious that a question is not simple search ."* and *" ...vice versa it is a problem of how to construct an algorithm that efficiently calculates $v_s$ [2, Section 3.3.4] "* are solved in excess of and the formulas at the level of elementary calculation are provided.

## §2. Spectral testing in the context of algebraic number theory

Spectral tests of R.Coveyou and R.MacPherson [4] for the checking of the arbitrary sequence $\{X_n\}$ with a maximum period of $N$ are built on the following ideas [1 and 2, Section 3.3.4].

For every positive integer $s$ on set $A_s \equiv \{t = (t_1, ..., t_s) : t_j = 0, ..., N-1 \ (j = 1, ..., s)\}$ the following characteristic of randomness is provided

$$f(t_1, ..., t_s) = \frac{1}{N}\sum_{k=0}^{N-1} \delta_{X_k}(t_1) \cdots \delta_{X_{k+(s-1)}}(t_s), \tag{2.1}$$

where here and below $\delta_M(x)$ is equal to 1 or 0 according to whether $x$ is a multiple number to $M$. This function is equal to the arithmetic mean (density) of occurrences of a particular combination $(t_1, ..., t_s)$ in terms of the $s$ consecutive sequence of $<X_n>: (t_1, ..., t_s) = (X_n, ..., X_{n+(s-1)})$.

In an ideal random sequence $<Y_n>$ with uniform distribution, all combinations $(t_1,...,t_s)$ must appear with the same frequency, so according to the corresponding (2.1) function $y(t_1,...,t_s) \equiv const$, where $const$ equal to $\frac{1}{N^s}$. Indeed, each combination $t = (t_1,...,t_s)$ in $A_s$ is exactly one, a number of the same such combinations is $N^s$, so the frequency in question, as the ratio of the number of occurrence of each to the number $t$ of all available opportunities $N^s$ is $\frac{1}{N^s}$.

Further, the direct final Fourier's transforms
$$\hat{f}(m_1,...,m_s) = \sum_{0 \le t_1,...,t_s < N} \exp\left(\frac{-2\pi i}{N}(m_1 t_1 + \cdots + m_s t_s)\right) f(t_1,...,t_s),$$
with the possibility of full unambiguous recovery of $f$ function by formulas of inverse finite Fourier's transforms
$$f(t_1,...,t_s) = \frac{1}{N^s} \sum_{\substack{m=(m_1,...,m_s) \in Z^s \\ 0 \le m_j < N (j=1,...,s)}} \exp\left(\frac{2\pi i}{N}(m_1 t_1 + \cdots + m_s t_s)\right) \hat{f}(m_1,...,m_s),$$
allow adequate transfer studies on Fourier, that means it deals with images of function $f$ with the use of this powerful analytical tool.

That is, the behavior of the Fourier coefficients of the function $f$, describing the distribution of the studied sequences $<X_n>$, is compared with the behavior of the Fourier coefficients of the function $y(t) \equiv \frac{1}{N^s}$, that describes perfectly the random sequence.

Calculations show (see. [1, Section 3.3.4A]) that for a linear congruential sequence
$$X_{n+1} = (aX_n + c) \bmod N \ (n = 0,1,2,...)$$
with a maximum period of $N$, Fourier's coefficients $f$ are equal to the following
$$\hat{f}(m_1,...,m_s) = \exp\left(\frac{-2\pi i c}{N} - \left(\frac{s(a,m) - s(1,m)}{a-1}\right)\right) \times \delta_N(m_1 + am_2 + \cdots + a^{s-1} m_s),$$
for (2.2) where $s(a,m) = m_1 + am_2 ... + a^{s-1} m_s$.

According to the definition of $\delta_N(x)$ from (2.2), it follows that $\hat{f}(m_1,...,m_s) \neq 0$ implies if and only if
$$m_1 + am_2 + \cdots + a^{s-1} m_s \equiv 0 (\bmod N), \tag{2.3}$$
and in all such cases when $|\hat{f}(m)| = 1$.

At the same time, the Fourier's coefficients $\hat{y}(m) \equiv \hat{y}(m_1,...,m_s)$ of functions $y(t)$, perfectly corresponding to the random sequence in the case of $m_j = 0 (\bmod N)$ for all $j = 1,...,s$, are equal to 1 and to 0 for all other $m$.

In conclusion, for $m \neq 0$ the values $\hat{f}(m)$ and $\hat{y}(m)$ are different only in values $\hat{f}(m) \neq 0$, at all points

$$\{m \in Z^s : m \neq 0, 0 \leq m_j < N(j = 1,...s), \hat{f}(m) \neq 0\} \equiv$$
$$\equiv \{m \in Z^s : m \neq 0, 0 \leq m_j < N(j = 1,...s), m_1 + am_2 + \cdots + a^{s-1}m_s \equiv 0 (\mod N)\} \quad (2.4)$$

Assumptions of R.Coveyou and R.MacPherson, constituting *spectral test*, are consisted in the fact that, the smallest Euclidean distance $\sqrt{m_1^2 + \cdots + m_s^2}$ of elements of the set (2.4) from zero $(0,...,0)$ should be taken as a measure of the deviation of a random sequence from the ideal possibility

$$v_s(a, N) = \min \sqrt{m_1^2 + \ldots + m_s^2}, \quad (2.5)$$

where the minimum is taken by all $s$ - set of integers $(m_1,\ldots,m_s) \neq (0,0,\ldots,0)$, which are the solutions of the congruence

$$m_1 + am_2 + \cdots + a^{s-1}m_s \equiv 0 (\mod N). \quad (2.6)$$

According to the [2, Section 3.3.4] (see. also [1, Section 3.3.4]), the definition (2.5)-(2.6) can be also given in the form of

$$v_s^2(a, N) = \inf \{(Nu_1 - au_2 - \ldots - a^{s-1}u_s)^2 + u_2^2 + \ldots + u_s^2 : (u_1,\ldots,u_s) \neq (0,\ldots,0)\}. \quad (2.7)$$

In this connection, we note that if in the spectral testing through the specific on $a$ and $N$ comparisons of (2.6), the "number" $m$ of non-zero Fourier's coefficients, measured in (2.5)-(2.7), are shifted as far as possible from zero, then in numerical integration the situation will be similar.

$a$ and $N$ should be such that the following inequality would be true $(\overline{m}_j = \max\{1; |m_j|\})$

$$\sum_{m_1,\ldots,m_s=-(N-1)}^{N-1} \frac{\delta_N(m_1 + am_2 + \cdots + a^{s-1}m_s)}{\overline{m}_1 \cdots \overline{m}_s} \leq \frac{\ln^\beta N}{N} (\beta > 0), \quad (2.8)$$

that can be satisfied if for all non-trivial solutions $m = (m_1,\ldots,m_s)$ of the comparison (2.6) the inequality $\overline{m}_1 \cdots \overline{m}_s \geq N$ is true (l. [5, p. 126-127]).

Thus, if in spectral comparison testing nontrivial solutions $m = (m_1,\ldots,m_s)$ of the comparison (2.6) should have large enough $\sqrt{m_1^2 + \cdots + m_s^2}$, which satisfies the definition of $v_s(a, N)$, the same value is required for the value $\overline{m}_1 \cdots \overline{m}_s$ in the numerical integration in the form of completion (2.8).

By this the following objectives can be achieved: in a spectral testing it is the building of good random number generators $\langle X_n(a, N) \rangle$, in numerical integration by $a$ and $N$ the nodes $\left(\left\{\frac{k}{N}\right\}, \left\{\frac{ak}{N}\right\}, \ldots, \left\{\frac{a^{s-1}k}{N}\right\}\right)(k=1,...,N)$ of good quadrature formula with equal weights are determined (here we must "avoid large Fourier's coefficients" of classes of functions with a dominant mixed derivative "numbers", in which $m$ of them form a so-called "hyperbolic crosses", for more details see in [5] and the current state, in [6]), where $\{x\}$ is a fractional part $x$.

In studies [7-12] the same problems of numerical integration were solved by the theory of divisors or through the same lattice (which is also done in this paper).

In this study, the basic idea is presenting of (2.7), which is for $a_j = a^{j-1} (j = 2,...,s)$ and taken from the following Theorem, that is another indication of the historical reliability of the mathematical nomenclature [13, p.429].

**Theorem C (K. Sherniyazov, [12]).** Let *a positive integer s and integers* $N \geq 2$, $a_1 = 1$, $a_2,...,a_s$ *are given, and let*

$$V_{N,a_2,...,a_s} = \begin{pmatrix} N & 0 & \cdots & 0 \\ -a_2 & 1 & \cdots & 0 \\ \vdots & \vdots & \ddots & \vdots \\ -a_s & 0 & \cdots & 1 \end{pmatrix}$$

*Then the following statements are true:*
*1. For any vector* $m = (m_1, m_2,...m_s) \in Z^s$, *satisfying the relation*

$$m_1 + a_2 m_2 + \cdots + a_s^{s-1} m_s \equiv 0 \pmod{N} \qquad (2.9)$$

*there is a vector,* $u = (u_1, u_2,...u_s) \in Z^s$ *such that* $m = u V_{N,a_2,...,a_s}$, *and for every m from* $Z^s$, *satisfying (2.9), such vector is unique.*
*2. Conversely, for any vector* $u = (u_1, u_2,...u_s) \in Z^s$ *the vector* $m = u V_{N,a_2,...,a_s}$ *is a solution of the congruence (2.9).*

Thus, all the evidence of this article can be regarded as closed in respect to studies [7-12] (of course, as it is said by modules of *Coveyou-MacPherson's method*).

## §3. The main results

We obtained the following results

**Theorem 1 (ST-2).** *Let there given an integer* $a \geq 5$ *and let* $(a-1)^2 = N$. *Then*

$$v_2(a, N; N = (a-1)^2) = \sqrt{1 + (a-2)^2} = (a-1) \times \sqrt{1 - 2\frac{a-2}{(a-1)^2}} =$$

$$= \sqrt{N - 2\sqrt{N} + 2} = \sqrt{N} \times \sqrt{1 - 2\frac{\sqrt{N} - 1}{N}} \approx a \approx \sqrt{N}$$

*and*

$$\mu_2(a; N; N = (a-1)^2) = \pi \left(1 - 2\frac{a-2}{(a-1)^2}\right) = \pi \left(1 - 2\frac{\sqrt{N} - 1}{N}\right).$$

**Theorem 2.** *Let the numbers* $s \geq 3$, $a \geq b_s + 1$ *and* $N > a$ *are related by* $(a-1)^s = N$. *Then the following relations will be held*

$$N^{\frac{1}{s}}\left(1 - (b_s - 1)N^{-\frac{1}{s}}\right) = a - b_s \leq v_s(a, N; N = (a-1)^s) \leq \sqrt{1 + a^2} = N^{\frac{1}{s}}\sqrt{1 + 2N^{-\frac{1}{s}} + 2N^{-\frac{2}{s}}}$$

*and*

$$\frac{\pi^{\frac{s}{2}}}{\left(\frac{s}{2}\right)!}\left(1-(b_s-1)N^{-\frac{1}{s}}\right)^s = \frac{\pi^{\frac{s}{2}}}{\left(\frac{s}{2}\right)!}\left(\frac{a-b_s}{a-1}\right)^s \leq \mu_s\left(a,N;N=(a-1)^s\right) \leq$$

$$\leq \frac{\pi^{\frac{s}{2}}}{\left(\frac{s}{2}\right)!} \cdot \frac{(1+a^2)^{\frac{s}{2}}}{(a-1)^s} = \frac{\pi^{\frac{s}{2}}}{\left(\frac{s}{2}\right)!} \cdot \left(1+2N^{-\frac{1}{s}}+2N^{-\frac{2}{s}}\right)^{\frac{s}{2}},$$

where $(-b_s)$ there is the biggest in absolute magnitude negative binomial coefficient in the expansion $(a-1)^s$ in powers $a$:

$$b_2=2, b_3=3, b_4=4, b_5=10, b_6=20, b_7=35, b_8=56, b_9=126, b_{10}=252, b_{11}=462, b_{12}=792,$$
$$b_{13}=1716, b_{14}=3432, b_{15}=6435,..., \text{etc.}$$

**Theorem 3** $(s=\tau=2, \lambda \geq 2)$. *Let the parameters $a$ and $N$ related by $(a-1)^2 = N\lambda$ equation, where the integer $\lambda \geq 2$. Then*

$$\sqrt{\frac{N}{\lambda}} \times \sqrt{1-2\frac{\sqrt{N\lambda}-1}{N\lambda}} = \frac{\sqrt{1+(a-2)^2}}{\lambda} \leq v_2\left(a,N;(a-1)^2=N\lambda, \lambda \geq 2\right) \overset{\lambda|(a-1)}{\leq} \frac{\sqrt{2}}{\lambda}(a-1) = \sqrt{\frac{2}{\lambda}} \cdot \sqrt{N},$$

*where is required $\lambda | (a-1) \Leftrightarrow (a-1) | N$ in upper bound of the condition.*

**Theorem 4** $(\tau > s=2, (a-1)^\tau = N)$. *Let there are given numbers $\tau > s = 2$, $a \geq b_\tau + 1$ and $N > a$ such that $N = (a-1)^\tau$. Then*

$$N^{\frac{1}{\tau}}\left(1-(b_\tau-1)N^{-\frac{1}{\tau}}\right) = a - b_\tau \leq v_2\left(a,N;(a-1)^\tau=N\right) \leq \sqrt{1+a^2} = N^{\frac{1}{\tau}}\sqrt{1+2N^{-\frac{1}{\tau}}+2N^{-\frac{2}{\tau}}},$$

$$\pi N^{\frac{2}{\tau}-1}\left(1-(b_\tau-1)N^{-\frac{1}{\tau}}\right)^2 = \pi\frac{(a-b_\tau)^2}{(a-1)^\tau} \leq \mu_2\left(a,N;(a-1)^\tau=N\right) \leq$$

$$\leq \pi\frac{(1+a^2)}{(a-1)^\tau} = \pi N^{\frac{2}{\tau}-1}\left(1+2N^{-\frac{1}{\tau}}+2N^{-\frac{2}{\tau}}\right).$$

**Theorem 5** $(\tau = s \geq 2, (a-1)^s = N\lambda, \lambda \geq 1)$. *Let the numbers $a \geq b_s+1, N > a, s \geq 2$ and $\lambda \geq 2$ are such that $(a-1)^s = N\lambda$. Thus*

$$\frac{N^{\frac{1}{s}}}{\lambda^{1-\frac{1}{s}}}\sqrt{\sum_{k=0}^{s-1}\left(\binom{s-1}{k}\right)^2} = \frac{a-1}{\lambda}\sqrt{\sum_{k=0}^{s-1}\left(\binom{s-1}{k}\right)^2} \overset{\lambda|(a-1)}{\geq} v_s\left(a,N;(a-1)^s=N\lambda\right) \geq$$

$$\geq \frac{a-b_s}{\lambda} = \frac{N^{\frac{1}{s}}}{\lambda^{1-\frac{1}{s}}}\left(1-(b_s-1)(N\lambda)^{-\frac{1}{s}}\right),$$

*where in the upper estimate assumes that the condition $\lambda|(a-1) \Leftrightarrow (a-1)^{\tau-1} | N$ holds.*

**Theorem 6** $(2 \leq s < \tau, (a-1)^\tau = N\lambda, 1 \leq \lambda \leq (a-1)^{\tau-s})$. *If $a, N, s, \tau$ and $\lambda$ are connected by the equality $(a-1)^\tau = N\lambda$ and inequalities $a \geq b_s+1, 2 \leq s < \tau, 1 \leq \lambda \leq (a-1)^{\tau-s}$, then*

$$(N\lambda)^{\frac{1}{\tau}}\left(1-(b_s-1)(N\lambda)^{-\frac{1}{\tau}}\right)=a-b_s\leq v_s\left(a,N;(a-1)^{\tau}=N\lambda\right)\leq\sqrt{1+a^2}=$$

$$=(N\lambda)^{\frac{1}{\tau}}\sqrt{1+2(N\lambda)^{-\frac{1}{\tau}}+2(N\lambda)^{-\frac{2}{\tau}}},$$

$$\frac{\pi^{\frac{s}{2}}}{\left(\frac{s}{2}\right)!}\frac{\left((N\lambda)^{\frac{1}{\tau}}\left(1-(b_s-1)(N\lambda)^{-\frac{1}{\tau}}\right)\right)^s}{N}=\frac{\pi^{\frac{s}{2}}}{\left(\frac{s}{2}\right)!}\frac{(a-b_s)^s\times\lambda}{(a-1)^{\tau}}\leq\mu_s\left(a,N;(a-1)^{\tau}=N\lambda\right)\leq$$

$$\leq\frac{\pi^{\frac{s}{2}}}{\left(\frac{s}{2}\right)!}\frac{(1+a^2)^{\frac{s}{2}}\times\lambda}{(a-1)^{\tau}}=\frac{\pi^{\frac{s}{2}}}{\left(\frac{s}{2}\right)!}\frac{\left((N\lambda)^{\frac{1}{\tau}}\sqrt{1+2(N\lambda)^{-\frac{1}{\tau}}+2(N\lambda)^{-\frac{2}{\tau}}}\right)^s}{N}.$$

**Theorem 7** $(2\leq\tau<s,(a-1)^{\tau}=N\lambda,\lambda\geq 1)$. *Let numbers $a,N>a,s,\tau$ and $\lambda$, that $(a-1)^{\tau}=N\lambda$, $2\leq\tau<s,\lambda\geq 1$ are given. Then*

$$v_s\left(a,N;(a-1)^{\tau}=N\lambda\right)\leq\sqrt{\sum_{k=0}^{\tau}\left(\binom{\tau}{k}\right)^2}.$$

### §4. Comments and conclusions

In this paragraph will specify some of the conditions in Theorems 1-7 and turn to a number of conclusions in the light of the results in them. We note immediately that there is no strict binding to the parameters of random number generators $a,N, 1\leq a<N$, that are acceptable for each application within allowed to change (in the sense that if the intended accuracy $1/10^9$ proved to be insufficient, then no one will be appointed to $1/(10^9+20110207)$ or $1/(10^9+22072013)$).

As with every pair of parameters $a$ and $N, a<N$, providing a maximum period of $N$ in (1.1)-(1.4) is uniquely associated number $\tau(a,N)\geq 2$ (see (1.3)), then returning to the conditions of $(a-1)^{\tau(a,N)}=N\lambda(a,N)$ in (1.3)-(1.4) by putting $\tau=\tau(a,N)$ and $\lambda=\lambda(a,N)$ in Theorems 1-7, and we will receive the corresponding results for the random number generator (1.1) with a maximum period of $N$, which we will not overwrite and will continue to use both proven theorems with preservation of their numbers.

In Theorems 3 and 5 it is assumed the fulfillment of condition $\lambda\,|\,(a-1)$ in the upper estimates, from which follows $N\!\not|\,(a-1)^{\tau(a,N)-1}$, i.e. the second condition in (1.3).

Fair conditions in these upper estimates show that the corresponding lower estimates, generally speaking, can not be improved in the sense of asymptotic order, so they can serve as a basis for the principal conclusions.

Thus, we will also assume the established ST-approval given at the beginning of this article.

Make preliminary conclusions from the results obtained here (see also §5) in the context of the concept of *"Only the true value $v_s$ determines the degree of randomness"*, according to what in [1-2], a large number of theoretical studies and computer searches was carried out, in which a random number generator i.e. multiplier $a$ and module $N$ fixed, and for $s = 2,...,6$ the $v_s(a,N), \lg v_s(a,N)$ and $\mu_s(a,N)$ calculated on them.

1°. In Theorem 1 were obtained the exact values $v_2$ and $\mu_2$ respectively for $s = 2$ and condition $N = (a-1)^2$. As for fixed $a$ and $N$ values $v_2(a,N) \geq v_3(a,N) \geq ...$ form a non-increasing sequence, then $v_2(a,N)$ - the largest of them.

The Theorem 1 (with additional confirmations in subsequent theorems 3 and 4) for $s = 2$ shows that with respect to random number generators with a maximum period, assessment [1, Section 3.3.4] and [2, Section 3.3.4] $v_2(a,N) \leq N^{\frac{1}{2}} \left(\frac{4}{3}\right)^{\frac{1}{4}}$ for all $N$ is overstated and that Indeed least upper bound $v_2(a,N)$ (greatest lower- too!) exactly is $N^{\frac{1}{2}}\sqrt{1 - 2\frac{\sqrt{N}-1}{N}}$ and for the coefficient the inequalities $\sqrt{1 - 2\frac{\sqrt{N}-1}{N}} < 1 < \left(\frac{4}{3}\right)^{\frac{1}{4}} = 1.07...$ is implemented.

Next, we turn to the results of a computer search of the [2, page 130.]: *the generators Lavaux and Janssens are located In lines 16 and 23; the parameters of these generators have been found on the computer to get a good multiplier in the sense of spectral criteria for which $\mu_2$ takes great importance:*

**Table 1.** Selected results of application of spectral criteria

| Line | $a$ | $N$ | $v_2^2$ | $\mu_2$ |
|------|-----------|----------|----------------|-------|
| 16   | 1664525   | $2^{32}$ | 4938916874     | 3.61  |
| 23   | 31167285  | $2^{48}$ | $3.2 \times 10^{14}$ | 3.60  |

Let us compare these data with the possibilities of Theorem 1. According to these multipliers $a$, we have:

Line 16: $a = 1664525, N = (a-1)^2 = 2770640146576$ and $v_2^2 = (a-2)^2 + 1 = 2770636817530$ against 4938916874 in Table 1.

Line 23: $a = 31167285, N = (a-1)^2 = 971399591936656$ and $v_2^2 = (a-2)^2 + 1 = 971399529602090 = 10^{14.987397...}$ against $3.2 \times 10^{14}$ in Table 1.

By the formula $\mu_2(a, N = (a-1)^2)$, according to ST-relations, may take arbitrarily close to its limit value - the $\pi$ number:

$$\left|\pi - \mu_2(a, N = (a-1)^2)\right| = 2\pi \frac{a-2}{(a-1)^2} < \frac{2\pi}{a-2}.$$

In particular, for multiplier $a = 1664525$ from the 16-th line we have
$\mu_2(a = 1664525, N = (a-1)^2 = 2770640146576) = \pi(1 - 0{,}000001201543984),$

against 3.61 Table 1, and, respectively, for $a = 31167285$ from 23 line:
$$\mu_2(a = 31167285, N = (a-1)^2 = 971399591936656) = \pi(1 - 2\frac{31167283}{971399591936656}) =$$
$$= \pi(1 - 0,000000064169850)$$
against 3.60 in Table 1.

2°. The answer to the crucial question in the spectral test: *"Since no one knows what the best achievable $v_s$ values"* [1, page 111.] obtained in Theorem 2:
$$v_s(a, N; (a-1)^s = N) \approx N^{\frac{1}{s}} \text{ and } \mu_s(a, N; (a-1)^s = N) \approx \pi^{\frac{s}{2}} / \left(\frac{s}{2}\right)!.$$

Moreover, as in the case $s = 2$, The Theorem 2 allows to specify $\gamma_N$ constants from Theorem B ([1, page Section 3.3.4] and [2, Section 3.3.4]) for $s \geq 3$:
if $N > \left(\frac{4}{\gamma^2(s)-1}\right)^s$, then $\overline{\gamma}_N(s) \equiv \sqrt{1 + 2 \cdot N^{-\frac{1}{s}} + 2N^{-\frac{2}{s}}} < \gamma(s)$, from what

$$\begin{cases} s = 3, & N \geq 3645 \Rightarrow \overline{\gamma}_N(3) < 2^{\frac{1}{6}}, \\ s = 4, & N \geq 8697 \Rightarrow \overline{\gamma}_N(4) < 2^{\frac{1}{4}}, \\ s = 5, & N \geq 28071 \Rightarrow \overline{\gamma}_N(5) < 2^{\frac{3}{10}}, \\ s = 6, & N \geq 47206 \Rightarrow \overline{\gamma}_N(6) < \left(\frac{64}{3}\right)^{\frac{1}{12}}, \\ s = 7, & N \geq 70730 \Rightarrow \overline{\gamma}_N(7) < 2^{\frac{3}{7}}, \\ s = 8, & N \geq 75628488 \Rightarrow \overline{\gamma}_N(8) < 2^{\frac{1}{2}}. \end{cases}$$

The same occurs for all $s \geq 9$ when the "coefficient" $\sqrt{1 + 2N^{-\frac{1}{s}} + 2N^{-\frac{2}{s}}}$ for $N^{\frac{1}{s}}$ in all these cases, in contrast to the case of $s = 2$ at the top infinitely close to 1 in unlimited quantities.

As part of the Theorem B after Theorems 1 and 2 in fact it would be possible to finish the search for "magic-magic" numbers with a maximum period, as they are actually received in the form of
$$v_s(a, N; N = (a-1)^s) \approx N^{\frac{1}{s}} = a - 1 \text{ and } \mu_s(a, N) \approx \pi^{\frac{s}{2}} / \left(\frac{s}{2}\right)!.$$

In the light of Theorem 1-2, the question of existence in the remaining cases 2°-4° of (1.8) of competitors to the case $N = (a-1)^s$ from the practical and theoretical points of view seem unpromising occupation as such it will be necessary to find $\gamma_N'$ and $\gamma_N''$, that
$$1 - (b_s - 1)N^{-\frac{1}{s}} = \overline{\overline{\gamma}}_N < \gamma'' \leq 1 \leq \gamma_N' < \overline{\gamma}_N = \sqrt{1 + 2N^{-\frac{1}{s}} + 2N^{-\frac{2}{s}}}.$$

However, in order to obtain a complete picture of the behavior of $s$-accuracy
$$v_s(a, N; (a-1)^{\tau(a,N)} = N\lambda(a, N))$$

these consideration (including subsequent non-trivial conclusions of the empirical data in [1-2]) continued in Theorems 3-7.

3°. Theorems 3 and 4 in the case $s=2$ about the influence of $\lambda$ in $(a-1)^s = N\lambda$ and a substitution of $s$ on $\tau$, $\tau > s$ at $(a-1)^s = N$ ($\tau < s$ impossible in $s=2$ case!).

From the Theorems 3 and 4 for $s=2$ follows that not $(a-1)^2$ increasing by $\lambda \geq 2$ in $(a-1)^2 = N\lambda$, not $\tau > s = 2$ increasing in $(a-1)^\tau = N$ not only will not give $v_2(a, N; (a-1)^2 = N) \approx \sqrt{N}$ increasing, even on the contrary, will lead to a decreasing for $\sqrt{\frac{2}{\lambda}}(\lambda \geq 2)$ and $N^{\frac{1}{2}-\frac{1}{\tau}}$ ($\tau > 2$) multipliers respectively.

Thus, in general, it seems to be in the $s=2$ case, when the greatest number of needs is provided, and even supplemented with an independence of $(X_n, X_{n+1})$ vectors, the best generator is $\langle X_n(a, N) \rangle$, which built on $a, N, a < N$ connected by the $(a-1)^2 = N$ equality.

The following Theorems are devoted to the general $s \geq 2$ case.

4°. Theorem 5, asymptotically exact for $\lambda | (a-1)$ gives accurate numerical information

$$v_s(a, N; (a-1)^s = N\lambda, 1 \leq \lambda < (a-1)^{s-1}) \approx \frac{1}{\lambda^{1-\frac{1}{s}}} \cdot N^{\frac{1}{s}},$$

and, simultaneously, shows the role of $\lambda$ in the $(a-1)^s = N\lambda$ ratio in accurate numerical data, namely deteriorates $v_s(a, N)$ growth for an accurate multiplier $1 / \lambda^{1-\frac{1}{s}}$.

The same happens with the $\mu_s(a, N)$ for $\lambda > 1$:

$$\mu_s(a, N; (a-1)^s = N\lambda, \lambda | (a-1)) \approx \frac{\pi^{\frac{s}{2}}}{\left(\frac{s}{2}\right)!} \cdot \frac{1}{\lambda^{s-1}},$$

in which asymptotically exact equality $\approx \pi^{\frac{s}{2}} / \left(\frac{s}{2}\right)!$ deteriorates by $1/\lambda^{s-1}$ multiplier for $\lambda = 1$.

5°. For $\lambda = 1$ theorem 6 gives an asymptotically exact solution of the problem of the construction of $< X_n(a, N) >$ generator with $v_2(a, N), \ldots, v_\tau(a, N)$ at the same time great for arbitrary $\tau > 2$. That is, for given $\tau > 2$ and $M > 0$ it is necessary to choose $a$ and $N$, $a < N$ such, that $(a-1)^\tau = N$ and $N^{\frac{1}{\tau}}(1-(b_\tau - 1) \cdot N^{-\frac{1}{\tau}}) > M$, and then for all $s$, $2 \leq s \leq \tau$, we have $v_s(a, N) \geq M$ that it follows from inequalities

$$N^{\frac{1}{\tau}}(1-(b_\tau-1)\cdot N^{-\frac{1}{\tau}}) \leq v_s(a, N; (a-1)^\tau = N) \leq N^{\frac{1}{\tau}}\sqrt{1+2N^{-\frac{1}{\tau}}+2N^{-\frac{2}{\tau}}} \quad (2 \leq s \leq \tau). \quad (4.1)$$

Final Theorem 7, in combination with previous, leads to the following recommendations.

6°. Conclusions on the $s$-dimensional accuracy of the random number generator with the maximum period. According to Theorems 1-7 for $\lambda = 1$ and given $\tau(a,N)$ for all $s, 2 \leq s \leq \tau(a,N)$ asymptotic equality

$$v_s(a, N; (a-1)^{\tau(a,N)} = N) \approx N^{\frac{1}{\tau(a,N)}}, \qquad (4.2)$$

is satisfied, which for $s > \tau(a,N)$ breaks to fall to

$$v_s(a, N; (a-1)^{\tau(a,N)} = N) \leq \sqrt{\sum_{k=0}^{\tau(a,N)} \left(\binom{\tau(a,N)}{k}\right)^2}. \qquad (4.3)$$

The confirmation of relations (4.1) - (4.3) is

**Table 1.** SELECTED RESULTS OF APPLICATION OF SPECTRAL CRITERIA

| Line | $a$ | $N$ | $v_2^2$ | $v_3^2$ | $v_4^2$ | $v_5^2$ | $v_6^2$ |
|---|---|---|---|---|---|---|---|
| 1 | 23 | $10^8+1$ | 530 | 530 | 530 | 530 | 447 |
| 2 | $2^7+1$ | $2^{35}$ | 16642 | 1664 | 16642 | 15602 | 252 |
| 3 | $2^{18}+1$ | $2^{35}$ | 34359738368 | 26 | 4 | 4 | 4 |

Here, in the second line $a-1=2^7$, $N=2^{35}$, so $\tau(a,N)=5$. In accordance with (4.1)-(4.2) the same values $v_2=v_3=v_4=16642$ with the peak $v_5^2=15602$ immediately, in accordance to (4.3), end up with falling to $v_6^2=252$.

For comparison we present $v_s(a=2^7+1, N=2^{35}, (a-1)^5=N)$ values calculated by Theorems 1 and 3:

| $a$ | $N$ | $\leq v_2^2$ | $\leq v_3^2$ | $\leq v_4^2$ | $\leq v_5^2$ | $\leq v_6^2$ |
|---|---|---|---|---|---|---|
| $2^7+1$ | $2^{35}$ | 16384 | 15876 | 15625 | 14161 | 252 |

Note that the 1-st row of Table 1, although $a-1=22=10^{1,342}, \tau=5,961\approx 6$ are not fully comply with the $N=10^8+1=(a-1)^\tau$ request, but for $s=2,3,4,5$ all $v_s^2$ are equal to 530 with little distortion for $s=6\approx\tau$.

The same applies to the 3rd line, for $\tau=\frac{35}{18}\approx 2$, and so $v_2=34359738368$, again in accordance with (4.3), following $v_s^2$ sharply fall:

$$v_3^2=6, v_4^2=v_5^2=v_6^2=4.$$

7°. Theorems 5 and 6 reveal an interesting role $\lambda>1$ in $v_s(a,N;(a-1)^\tau=N\lambda)$: if for $\tau=s$ the parameter $\lambda>1$ leads to a $N^{\frac{1}{\tau}}$ decrease by a $1/\lambda^{1-\frac{1}{s}}$ multiplier, then for $\tau>s$, $1\leq\lambda\leq(a-1)^{\tau-s}$ to increasing in $\lambda^{\frac{1}{\tau}}$ times.

In the result - a general recommendation on $\lambda\geq 1$: if $v_s(a,N), \mu_s(a,N)$ are necessary for one $s\geq 2$ to choose the highest, then, $a$ and $N$, $a<N$ should be connected by the equations $\lambda=1$ and $(a-1)^s=N$, if the same is required for $s, s+1,\ldots s+l,$ then the selection $\tau$ and $\lambda$ in $(a-1)^\tau=N\lambda$ should be in the ratios of $\tau>s+l$ and $1\leq\lambda\leq(a-1)^{\tau-(s+l)}$.

We combine in one general all of the proposed private recommendations on the building of the best possible generators.

8°. Rule of constructing of $\langle X_n(a,N) \rangle$ random numbers generator with the maximum period by Theorems 1-7.

1. If for given $s, s \geq 2$ it is required to construct the generator with the greatest $v_s(a,N)$ and $\mu_s(a,N)$, then $a$ and $N$, $a < N$ must be connected by the equation $(a-1)^s = N$, and then

$$v_2(a,N; N=(a-1)^2) = \sqrt{1+(a-2)^2} = (a-1) \times \sqrt{1 - 2\frac{a-2}{(a-1)^2}} =$$

$$= \sqrt{N - 2\sqrt{N} + 2} = \sqrt{N} \times \sqrt{1 - 2\frac{\sqrt{N}-1}{N}} \approx a \approx \sqrt{N},$$

$$\mu_2(a; N; N=(a-1)^2) = \pi\left(1 - 2\frac{a-2}{(a-1)^2}\right) = \pi\left(1 - 2\frac{\sqrt{N}-1}{N}\right),$$

and further for all $s, s \geq 3$

$$N^{\frac{1}{s}}\left(1 - (b_s - 1)N^{-\frac{1}{s}}\right) = a - b_s \leq v_s(a,N; N=(a-1)^s) \leq \sqrt{1+a^2} = N^{\frac{1}{s}}\sqrt{1 + 2N^{-\frac{1}{s}} + 2N^{-\frac{2}{s}}},$$

$$\frac{\pi^{\frac{s}{2}}}{\left(\frac{s}{2}\right)!}\left(1 - (b_s - 1)N^{-\frac{1}{s}}\right)^s = \frac{\pi^{\frac{s}{2}}}{\left(\frac{s}{2}\right)!}\left(\frac{a-b_s}{a-1}\right)^s \leq \mu_s(a,N; N=(a-1)^s) \leq$$

$$\leq \frac{\pi^{\frac{s}{2}}}{\left(\frac{s}{2}\right)!} \cdot \frac{(1+a^2)^{\frac{s}{2}}}{(a-1)^s} = \frac{\pi^{\frac{s}{2}}}{\left(\frac{s}{2}\right)!} \cdot \left(1 + 2N^{-\frac{1}{s}} + 2N^{-\frac{2}{s}}\right)^{\frac{s}{2}}.$$

Note that for a fixed $s$ $(a-1)^s = N$ condition provided the opportunity to be an arbitrarily large to a multiplier $a = d \times p_1^{r_1} \ldots p_t^{r_t} + 1$ and the maximum period $N = d^s \times p_1^{sr_1} \ldots p_t^{sr_t}$.

2. If for given $\tau \geq 3$ it is necessary to build generator "with large (guaranteed not less) than $(N\lambda)^{\frac{1}{\tau+l}} \approx a - b_{\tau+l}$" $v_2(a,N), \ldots, v_\tau(a,N)$, then choosing $a, N, l \geq 1$, $1 \leq \lambda \leq (a-1)^l = (N\lambda)^{\frac{l}{\tau+l}}$ those that $(a-1)^{\tau+l} = N\lambda$, then we get

$$(N\lambda)^{\frac{1}{\tau+l}} \times \left(1 - (b_{\tau+l} - 1)(N\lambda)^{-\frac{1}{\tau+l}}\right) = a - b_{\tau+l} \leq v_s(a,N;(a-1)^{\tau+l} = N\lambda) \leq$$

$$\leq \sqrt{a^2 + 1} = (N\lambda)^{\frac{1}{\tau+l}}\sqrt{1 + 2(N\lambda)^{-\frac{1}{\tau+l}} + 2(N\lambda)^{-\frac{2}{\tau+l}}}$$

for all $2 \leq s \leq \tau$.

3. If for given $\tau \geq 2$ it is necessary to build generator "with large (guaranteed not less) than $\approx \dfrac{\pi^{\frac{s}{2}}}{\left(\dfrac{s}{2}\right)!} \cdot \dfrac{\lambda^{\frac{1}{\tau+l}}}{N^{1-\frac{1}{\tau+l}}}) \gg \mu_2(a,N),...,\mu_\tau(a,N)$, then choosing $a, N, l$ and $\lambda$ those that $(a-1)^{\tau+l} = N\lambda, 1 \leq \lambda \leq (a-1)^l = (N\lambda)^{\frac{l}{\tau+l}}$, then we get

$$\dfrac{\pi^{\frac{s}{2}}}{\left(\dfrac{s}{2}\right)!} \cdot \dfrac{(N\lambda)^{\frac{s}{\tau+l}}}{N}\left(1-(b_s-1)(N\lambda)^{-\frac{1}{\tau+l}}\right)^s \leq$$

$$\leq \mu_s\left(a, N; (a-1)^{\tau+l} = N\lambda\right) = \dfrac{\pi^{\frac{s}{2}}}{\left(\dfrac{s}{2}\right)!} \cdot \dfrac{v_s^s\left(a, N; (a-1)^{\tau+l} = N\lambda\right)}{N} \leq$$

$$\leq \dfrac{\pi^{\frac{s}{2}}}{\left(\dfrac{s}{2}\right)!} \cdot \dfrac{(N\lambda)^{\frac{s}{\tau+l}}}{N}\left(1+2(N\lambda)^{-\frac{1}{\tau+l}}+2(N\lambda)^{-\frac{2}{\tau+l}}\right)^{\frac{s}{2}}.$$

for all $2 \leq s \leq \tau$.

Of course, this calculation by $\lambda^{\frac{s}{\tau+l}}, 1 \leq \lambda \leq (a-1)^l$ multiplier. In this case, a further increase of $\lambda$ in general does not lead to an $\mu_s\left(a, N; (a-1)^{\tau+l} = N\lambda\right)$ increase.

9°. Thus, the asymptotic formula *ST* is theoretically almost explicit asymptotic estimates provides complete freedom of choice of $N, s$ and $a$ with a very optimal indexes. Therefore, the efficiency of application is transferred to the technical capabilities of constantly improved computers.

Here we must not to lose from view that, in contrast to the unlimited theoretical and computational capabilities *ST* relations, then arbitrariness is not quite large. Indeed, as noted in [3, p.34], *"If you consider that the lifetime of the universe is approximately 10 billion = $10^{10}$ years $< 10^{18}$ seconds, it is clear that none of the most fantastic calculation speeds will not provide the required accuracy"*, in our case $N = 10^{18+k}$, where $k$ - number of computer operations in one second, which currently stands as 33.8*1015 operations per second (according to http://top500.org/blog/lists/2015/11/press-release/).

10°. In Theorems 1-7 for a fixed $s, s \geq 2$ asymptotic behavior $v_s(a,N)$ determined depending on relationships $\tau = s \geq 2, \tau > s$ and $2 \leq \tau < s$, where $\tau$ is a potential of the $\langle X_n(a,N) \rangle$ generator binding multiplier $a$ and a maximum period $N$ of $(a-1)^\tau = N\lambda$ equality.

However, in applications as a result of Theorems 1-7, as it was shown in the $3° - 9°$ roles of $\tau$ and $s$ are changing - the random number generator $\langle X_n(a,N) \rangle_{n=0}^{N-1}$ is defined by $\tau, a$, its spectral testing is performed according to the values of the quantities $v_s(a, N; (a-1)^\tau = N\lambda)$ for a variable $s \geq 2$, asymptotically exact orders, which in advance are found in the same theorems 1 -7.

11°. Theorems 1-7 lead to the general conclusion: $a, a \approx v_s(a, N)$ multiplier responsible for $s$-accuracy and potential $\tau$ - of the value of the maximum period $N\lambda = (a-1)^\tau$ and for effective limits $s$, $2 \le s \le \tau$.

12°. In open applications for the given multiplier $a - 1 = d \times p_1^{r_1} \cdots p_t^{r_t}$ and the potential $\tau \ge 2$ we can limit ourselves with the $N = (a-1)^\tau = p_1^{\pi_1} \cdots p_t^{\pi_t}$ case, while providing almost unlimited confidentiality maximum period of $N = p_1^{\aleph_1} \cdots p_t^{\aleph_t}$, combined with $\lambda = d \times p_1^{\pi_1 - \aleph_1} \cdots p_t^{\pi_t - \aleph_t}$ shall be selected within the ST-limits with

$$2 \le s \le \tau, \quad a - b_\tau \le v_s\left(a, N; (a-1)^\tau = N\lambda\right) \le \sqrt{a^2 + 1}$$

Estimates.

13°. In applications it is necessary not to lose sight of the difference between the asymptotically exact formulas and specific calculations on them.

In Theorems 1-7 assessments have the form of $N^{\frac{1}{s}} \underline{\underline{\gamma}}_N \le v_s(a, N) \le N^{\frac{1}{s}} \overline{\gamma}_N$, that carry asymptotic character with a clear record of the values $\overline{\gamma}_N$ and $\overline{\overline{\gamma}}_N$.

Absolutely accurate for an unlimited increase of $N$ in the specific calculations with typically low $N$ values, errors can be appreciable.

In this connection, we note that in the ST-method $\overline{\overline{\gamma}}_N = \left(1 - (b_\tau - 1)N^{-\frac{1}{\tau}}\right)$ and $\overline{\gamma}_N = \sqrt{1 + 2N^{-\frac{1}{\tau}} + 2N^{-\frac{2}{\tau}}}$ multipliers for $N^{\frac{1}{\tau}}$ are only guaranteeing closeness to the best value $s$-accuracy $v_s \equiv v_s(a, N; (a-1)^\tau = N)$ random $\langle X_n(a, N = (a-1)^\tau)\rangle_{n=0}^{N-1}$ consequence, so although in applications written $... \le v_s$, the reality may be far exceeding those guarantees.

14°. Selecting the maximum period $N = (a-1)^\tau = p_1^{\aleph_1} \cdots p_t^{\aleph_t}$ is arbitrary due to the arbitrary choice of simple $1 < p_1 < p_2 < ... < p_t$ and positive $\aleph_1, ..., \aleph_t$, arbitrary choice of multiplier $a = d \times p_1^{r_1} \cdots p_t^{r_t} + 1$ with an arbitrary positive $r_1, ..., r_t$ and arbitrary choice of potential $\tau \ge 2$, but related with $\aleph_j \le r_j \tau$ inequalities for all $j = 1, ..., s$, and $(\tau - 1)r_{j_0} < \aleph_{j_0}$ for some $j_0$, that is consistent with the requirements in [2, page 118.]:

*"... when the applications need a random number generator, which provides reception of a consequence, which is very close to random, simple congruent generators for this purpose are not suitable. Instead, they need to use a generator with a long period, even if in fact it is necessary to generate only a small part of the period ".*

That is, for a fixed multiplier $a = d \times p_1^{r_1} \cdots p_t^{r_t} + 1$ by $\tau$ potential increase it is possible to achieve an arbitrarily large value of the maximum period $N = p_1^{\pi_1} \cdots p_t^{\pi_t}$. And not just "big", but at the same time ensuring the quality of controlled random sequence $\langle X_n(a, (a-1)^\tau = N\lambda)\rangle_{n=0}^{N-1}$ in the form of independence of $s$-dimensional vectors $(X_n, X_{n+1}, ..., X_{n+(s-1)})$ for all $2 \le s \le \tau$ with bilateral estimates of the asymptotic type $v_s(a,N;(a-1)^\tau = N\lambda) \approx (N\lambda)^{\frac{1}{\tau}}$ and $\mu_s(a,N;(a-1)^\tau = N\lambda) \approx \pi^{s/2}/(s/2)!N \times (N\lambda)^{s/\tau}$.
Thus here indicators are extremely accurate.

15°. Selecting an *a* multiplier allowing for the specials of computer on which the calculations are carried out, are discussed in [2, 3.2.1.3. Potential.] *"In the previous section it was shown that the maximum period can be achieved when $a-1$ multiple to each prime divisor $N$, and $a-1$ also must be a multiple of 4 if $N$ multiple of 4. If b - machine base radix ($b=2$ for binary computer and $b=10$ - for decimal computer), $N$ - the word length in the computer $b^t$, then the multiplier $a = b^k + 1, 2 \leq k < t$ satisfies these conditions. By Theorem A it can be taken $c=1$. Recursive relation now has the form*
$$X_{n+1} = ((b^k + 1)X_n + 1) \bmod b^t$$
*and this equation means that the multiplication can be avoided; simply move and summation."*

Then, according to the ST relations, for $N = b^t$ we get $(a-1)^\tau = b^{k\tau} = b^t = N$ for all $\tau$ that share $t$, with two-way assessment
$$a - b_\tau \leq v_s(a,N;(a-1)^\tau = N) \leq \sqrt{1+a^2}$$
for all $2 \leq s \leq \tau$.

Further, in [2, Section 3.2.1.3] reported "For example, let $a = B^2 + 1$, where B - the size of the MIX computer bytes. Program LDAX; 2 SLA; ADDX; INCA 1 can be used instead of the program referred to in section 3.2.1.1, and the program execution time is reduced from 16u up to 7u.

*For this reason, multipliers having the $a = b^k + 1$ form widely discussed in the literature. They actually recommended by many authors. However, the first few years of experimenting with this method have convincingly shown that the multipliers that have a simple form $a = b^k + 1$ should be avoided. The generated random number is simply not enough."*

We can assume that accurate, instead of experimental, ST-ratios in each case lead to the right (adequate) conclusions.

16°. Random selection of $\tau$ potential support for all $s, 2 \leq s \leq \tau$ the independence of a $s$-sequence of values of a linear congruent sequence $\{X_n\}$ provided bilaterally asymptotically precise $\approx N^{\frac{1}{\tau}}$ estimate on a ST formula that removes all restrictions related to the *"available for [1-2]"* cases $s = 2,3,4,5,6$ (see. [1 and 2, Section 3.3 .4]).

17°. Checking the fulfillment of conditions (1.3)-(1.4) for a given *a* and *N* is not a trivial task. For example, in [1, Section 3.2.1.3] and [2, Section 3.2.1.3] introduced specific definitions, in which the smallest solution $\tau$ of the congruence $a^x \equiv 1 (\bmod N)$ called the order multiplier *a* in module *N*, and any such value *a* that has the *maximum* possible order in module *N* - *a primitive element* in module *N*.

Apparently, Theorems 1-7, as shown in 5°-16°, eliminate the need for such studies, it is only in cases dictated by other considerations.

18°. It is possible that a combination of creation of new methods by known generators (see [1, Section 3.2.1.3] and [2, Section 3.3.4]) using the constructed herein may be promising.

If with $s$-accuracy random number generator $v_s \approx N^{\frac{1}{\tau(a,N)}}$ for $2 \leq s \leq \tau(a,N)$ a more or less well, then, as shown by ST-asymptotically precise equality for fixed $a$ and $N$ where $N$ is a maximum period of a random number generator, generally speaking, at the same time ensuring high $\mu_s(a,N)$ for the same $2 \leq s \leq \tau(a,N)$, to the possible extent require separate computing in the context of Theorem 6 (see paragraph 3 of 8°).

There is a view of the following:

19°. Measure of the effectiveness $\mu_s(a,N)$ of the multiplier $a$ for a given module $N$, cited above, was introduced as a "*relatively independent of $N$*" rule for determining the quality of the random number generator.

Since, as shown by Theorems 1-7, the upper bound $\leq \gamma_N N^{\frac{1}{s}}$ for the $v_s(a,N)$ is achievable, it is natural to define $\mu_s(a,N)$ as the ratio $\frac{v_s^s(a,N)}{N}$, without the multiplier $\pi^{\frac{s}{2}} / \left(\frac{s}{2}\right)!$.

This proposal is based on the fact, that this multiplier has emerged as the volume of the ellipsoid in $s$-dimensional space, but that, as shown by a study carried out here has not directly relation to the subject under discussion.

However, we continue in the accepted definitions and notations of comparison (for $\lambda = 1$) with known [1-2].

In full accord with $v_s^s\left(a, N; N = (a-1)^{\tau(a,N)}\right) \approx N^{\frac{s}{\tau(a,N)}}$ the measure of efficiency $a$ and $N$ increases with an increase of $s$ from 2 up to $\tau(a,N)$:

$$\mu_s\left(a, N; N = (a-1)^{\tau(a,N)}\right) = \frac{\pi^{\frac{s}{2}} v_s^s\left(a, N; N = (a-1)^{\tau(a,N)}\right)}{\left(\frac{s}{2}\right)! N} \approx$$

$$\approx \frac{\pi^{\frac{s}{2}}}{\left(\frac{s}{2}\right)!} \cdot \frac{1}{N^{1-\frac{s}{\tau(a,N)}}} \quad (2 \leq s \leq \tau(a,N)),$$

(4.4)

and increases to a maximum $\dfrac{\pi^{\frac{\tau(a,N)}{2}}}{\left(\frac{\tau(a,N)}{2}\right)!}$ and then, when passing $s$ through $\tau(a,N)$ it becomes "*badly small*"

$$\mu_s(a,N) \le \frac{\pi^{\frac{s}{2}} \cdot \left(\sum_{k=0}^{\tau}\left(\binom{\tau(a,N)}{k}\right)^2\right)^{\frac{s}{2}}}{\left(\frac{s}{2}\right)!} \cdot \frac{1}{N} \quad (s > \tau(a,N)). \tag{4.5}$$

These conclusions are also supported by the experimental data [2, Section 3.3.4]

| Line | $a$ | $N$ | $\mu_2$ | $\mu_3$ | $\mu_4$ | $\mu_5$ | $\mu_6$ | |
|---|---|---|---|---|---|---|---|---|
| 2 | $2^7+1$ | $2^{35}$ | $2e^6$ | $3e^4$ | 0.04 | 4.66 | $2e^3$ | $e=\frac{1}{10}$ |

Indeed, in accordance with (4.4)-(4.5) for $\tau=5$ occurs rise of $\mu_2 = \frac{2}{10^6} < \mu_3 = \frac{3}{10^4} < \mu_4 = \frac{4}{10^2}$ with a peak $\mu_5 = 4.66$ and sharp fall of $\mu_6 = \frac{2}{10^3}$.

Numerical data obtained from the estimates (4.4)-(4.5) confirm or deny (then, apparently, it is necessary to recalculate) these figures:

| $a=2^7$ | $N=2^{35}$ | $1.47471e^6 \le$ $\le \mu_2 \le$ $\le 1.52161e^6$ | $2.4386551e^4 \le$ $\le \mu_3 \le$ $\le 2.6172608e^4$ | $3.506387856e^2 \le$ $\le \mu_4 \le$ $\le 3.977690029e^2$ | $3.65580844757 \le$ $\le \mu_5 \le$ $\le 5.47346593181$ | $\mu_6 \le 2.40685602e^3$ |
|---|---|---|---|---|---|---|

20°. In all tasks that are required to construct a uniformly distributed on $[0,1]$ sequence

$$U_n = \frac{X_n}{N} \quad (n=0,...,N-1)$$

ensuring $s$-random in the sense of $s$ independence of consecutive $(U_n, U_{n+1}, ..., U_{n+(s-1)})$ numbers on the full period of $N$ members and if possible to restrict with accuracy $\frac{1}{v_s(a, N=(a-1)^s)}$ when the first $\lg v_s(a, N=(a-1)^s)$ bits in the binary representation of the numbers can be considered random, then it is need to take $\langle X_n(a, N=(a-1)^s, c=1, X_0=0)\rangle_{n=0}^{N-1}$ as $\langle X_n\rangle$.

Checking the adequacy of the choice of "magic-magic numbers $N, a, c, X_0$" and $s$ to achieve this goal, in addition to the specific requirements of the task assigned, should also include a test of the independence of the [1-2] methods and tested for uniformity of distribution, for example, by calculating the difference module of frequency entering $\langle U_n\rangle$ in the arbitrary segment $[\alpha,\beta]$ and its length $\beta-\alpha$.

Lets start computational experiments with random number generator, built by $a=26, \tau=2$ and $N=625$:

$$X_{n+1} = (X_n \cdot 26 + 1) \bmod 625 \, (n=0,...,624; X_0=0).$$

The resulting sequence of random numbers $\langle U_n\rangle_{n=0}^{N-1}$ taken out in Table 2:

**Table 2**. Random numbers.

0,0016; 0,0432; 0,1248; 0,2464; 0,408; 0,6096; 0,8512; 0,1328; 0,4544; 0,816;
0,2176; 0,6592; 0,1408; 0,6624; 0,224; 0,8256; 0,4672; 0,1488; 0,8704; 0,632;
0,4336; 0,2752; 0,1568; 0,0784; 0,04; 0,0416; 0,0832; 0,1648; 0,2864; 0,448;
0,6496; 0,8912; 0,1728; 0,4944; 0,856; 0,2576; 0,6992; 0,1808; 0,7024; 0,264;
0,8656; 0,5072; 0,1888; 0,9104; 0,672; 0,4736; 0,3152; 0,1968; 0,1184; 0,08;
0,0816; 0,1232; 0,2048; 0,3264; 0,488; 0,6896; 0,9312; 0,2128; 0,5344; 0,896;
0,2976; 0,7392; 0,2208; 0,7424; 0,304; 0,9056; 0,5472; 0,2288; 0,9504; 0,712;
0,5136; 0,3552; 0,2368; 0,1584; 0,12; 0,1216; 0,1632; 0,2448; 0,3664; 0,528;
0,7296; 0,9712; 0,2528; 0,5744; 0,936; 0,3376; 0,7792; 0,2608; 0,7824; 0,344;
0,9456; 0,5872; 0,2688; 0,9904; 0,752; 0,5536; 0,3952; 0,2768; 0,1984; 0,16;
0,1616; 0,2032; 0,2848; 0,4064; 0,568; 0,7696; 0,0112; 0,2928; 0,6144; 0,976;
0,3776; 0,8192; 0,3008; 0,8224; 0,384; 0,9856; 0,6272; 0,3088; 0,0304; 0,792;
0,5936; 0,4352; 0,3168; 0,2384; 0,2; 0,2016; 0,2432; 0,3248; 0,4464; 0,608;
0,8096; 0,0512; 0,3328; 0,6544; 0,016; 0,4176; 0,8592; 0,3408; 0,8624; 0,424;
0,0256; 0,6672; 0,3488; 0,0704; 0,832; 0,6336; 0,4752; 0,3568; 0,2784; 0,24;
0,2416; 0,2832; 0,3648; 0,4864; 0,648; 0,8496; 0,0912; 0,3728; 0,6944; 0,056;
0,4576; 0,8992; 0,3808; 0,9024; 0,464; 0,0656; 0,7072; 0,3888; 0,1104; 0,872;
0,6736; 0,5152; 0,3968; 0,3184; 0,28; 0,2816; 0,3232; 0,4048; 0,5264; 0,688;
0,8896; 0,1312; 0,4128; 0,7344; 0,096; 0,4976; 0,9392; 0,4208; 0,9424; 0,504;
0,1056; 0,7472; 0,4288; 0,1504; 0,912; 0,7136; 0,5552; 0,4368; 0,3584; 0,32;
0,3216; 0,3632; 0,4448; 0,5664; 0,728; 0,9296; 0,1712; 0,4528; 0,7744; 0,136;
0,5376; 0,9792; 0,4608; 0,9824; 0,544; 0,1456; 0,7872; 0,4688; 0,1904; 0,952;
0,7536; 0,5952; 0,4768; 0,3984; 0,36; 0,3616; 0,4032; 0,4848; 0,6064; 0,768;
0,9696; 0,2112; 0,4928; 0,8144; 0,176; 0,5776; 0,0192; 0,5008; 0,0224; 0,584;
0,1856; 0,8272; 0,5088; 0,2304; 0,992; 0,7936; 0,6352; 0,5168; 0,4384; 0,4;
0,4016; 0,4432; 0,5248; 0,6464; 0,808; 0,0096; 0,2512; 0,5328; 0,8544; 0,216;
0,6176; 0,0592; 0,5408; 0,0624; 0,624; 0,2256; 0,8672; 0,5488; 0,2704; 0,032;
0,8336; 0,6752; 0,5568; 0,4784; 0,44; 0,4416; 0,4832; 0,5648; 0,6864; 0,848;
0,0496; 0,2912; 0,5728; 0,8944; 0,256; 0,6576; 0,0992; 0,5808; 0,1024; 0,664;
0,2656; 0,9072; 0,5888; 0,3104; 0,072; 0,8736; 0,7152; 0,5968; 0,5184; 0,48;
0,4816; 0,5232; 0,6048; 0,7264; 0,888; 0,0896; 0,3312; 0,6128; 0,9344; 0,296;
0,6976; 0,1392; 0,6208; 0,1424; 0,704; 0,3056; 0,9472; 0,6288; 0,3504; 0,112;
0,9136; 0,7552; 0,6368; 0,5584; 0,52; 0,5216; 0,5632; 0,6448; 0,7664; 0,928;
0,1296; 0,3712; 0,6528; 0,9744; 0,336; 0,7376; 0,1792; 0,6608; 0,1824; 0,744;
0,3456; 0,9872; 0,6688; 0,3904; 0,152; 0,9536; 0,7952; 0,6768; 0,5984; 0,56;
0,5616; 0,6032; 0,6848; 0,8064; 0,968; 0,1696; 0,4112; 0,6928; 0,0144; 0,376;
0,7776; 0,2192; 0,7008; 0,2224; 0,784; 0,3856; 0,0272; 0,7088; 0,4304; 0,192;
0,9936; 0,8352; 0,7168; 0,6384; 0,6; 0,6016; 0,6432; 0,7248; 0,8464; 0,008;
0,2096; 0,4512; 0,7328; 0,0544; 0,416; 0,8176; 0,2592; 0,7408; 0,2624; 0,824;
0,4256; 0,0672; 0,7488; 0,4704; 0,232; 0,0336; 0,8752; 0,7568; 0,6784; 0,64;
0,6416; 0,6832; 0,7648; 0,8864; 0,048; 0,2496; 0,4912; 0,7728; 0,0944; 0,456;
0,8576; 0,2992; 0,7808; 0,3024; 0,864; 0,4656; 0,1072; 0,7888; 0,5104; 0,272;
0,0736; 0,9152; 0,7968; 0,7184; 0,68; 0,6816; 0,7232; 0,8048; 0,9264; 0,088;
0,2896; 0,5312; 0,8128; 0,1344; 0,496; 0,8976; 0,3392; 0,8208; 0,3424; 0,904;

0,5056; 0,1472; 0,8288; 0,5504; 0,312; 0,1136; 0,9552; 0,8368; 0,7584; 0,72;
0,7216; 0,7632; 0,8448; 0,9664; 0,128; 0,3296; 0,5712; 0,8528; 0,1744; 0,536;
0,9376; 0,3792; 0,8608; 0,3824; 0,944; 0,5456; 0,1872; 0,8688; 0,5904; 0,352;
0,1536; 0,9952; 0,8768; 0,7984; 0,76; 0,7616; 0,8032; 0,8848; 0,0064; 0,168;
0,3696; 0,6112; 0,8928; 0,2144; 0,576; 0,9776; 0,4192; 0,9008; 0,4224; 0,984;
0,5856; 0,2272; 0,9088; 0,6304; 0,392; 0,1936; 0,0352; 0,9168; 0,8384; 0,8;
0,8016; 0,8432; 0,9248; 0,0464; 0,208; 0,4096; 0,6512; 0,9328; 0,2544; 0,616;
0,0176; 0,4592; 0,9408; 0,4624; 0,024; 0,6256; 0,2672; 0,9488; 0,6704; 0,432;
0,2336; 0,0752; 0,9568; 0,8784; 0,84; 0,8416; 0,8832; 0,9648; 0,0864; 0,248;
0,4496; 0,6912; 0,9728; 0,2944; 0,656; 0,0576; 0,4992; 0,9808; 0,5024; 0,064;
0,6656; 0,3072; 0,9888; 0,7104; 0,472; 0,2736; 0,1152; 0,9968; 0,9184; 0,88;
0,8816; 0,9232; 0,0048; 0,1264; 0,288; 0,4896; 0,7312; 0,0128; 0,3344; 0,696;
0,0976; 0,5392; 0,0208; 0,5424; 0,104; 0,7056; 0,3472; 0,0288; 0,7504; 0,512;
0,3136; 0,1552; 0,0368; 0,9584; 0,92; 0,9216; 0,9632; 0,0448; 0,1664; 0,328;
0,5296; 0,7712; 0,0528; 0,3744; 0,736; 0,1376; 0,5792; 0,0608; 0,5824; 0,144;
0,7456; 0,3872; 0,0688; 0,7904; 0,552; 0,3536; 0,1952; 0,0768; 0,9984; 0,96;
0,9616; 0,0032; 0,0848; 0,2064; 0,368; 0,5696; 0,8112; 0,0928; 0,4144; 0,776;
0,1776; 0,6192; 0,1008; 0,6224; 0,184; 0,7856; 0,4272; 0,1088; 0,8304; 0,592;
0,3936; 0,2352; 0,1168; 0,0384; 0

For this segment $[\alpha, \beta] (0 \leq \alpha < \beta \leq 1)$ we denote $m \equiv m(\alpha, \beta)$ by the number of random numbers $U_n$ belonging to a segment $[\alpha, \beta]$ and, accordingly, the frequency $\frac{m}{N} = \frac{m(\alpha, \beta)}{N}$ of falling $\langle U_n \rangle$ into an segment $[\alpha, \beta]$ value is then taken as a measure of the uniform distribution - the smaller, the better. We will be guided by the following examples sections select: [*year, month, day, father's or mother's birth; year, month, day of the student's birth*], $\left[\frac{1}{\pi^2}, 1 - \frac{1}{e}\right]$ and $[0.2; 0.9]$.

**Table 3.** The test results of random numbers in Table 2 for the uniform distribution.

|  | α | β | m | m/N | β-α | Δ |
|---|---|---|---|---|---|---|
| $\langle X_n(a = 26; N = 625)\rangle_{n=0}^{624}$ | 0,580815 | 0,850411 | 168 | 0,2688 | 0,2696 | 0,000796 |
|  | 1/π² | 1-1/e | 332 | 0,5312 | 0,530799 | 0,0004 |
|  | 0,2 | 0,9 | 437 | 0,6992 | 0,7 | 0,0008 |

The results are encouraging – Δ with three zeros after the decimal point with a relatively small $N = 625$.

Thus, the particular generator
$$X_{n+1} = (X_n \cdot 26 + 1) \bmod 625 \, (n = 0, \ldots, 624; X_0 = 0)$$
stood the test in the uniform distribution.

Regarding the tests for independence in advance, you can assume that any of random numbers $ST$-generator, including of course, and reporting, has passed the test criterion series **[2, * F. Connection with the criterion of the series]**.

## §5. Examples of the random number generator with a maximum period in the context of well-known and popular

The title of this section is devoted to the Table 1 in [2, Section 3.3.4], where in 29 examples with detailed commentary the results of theoretical studies and computer searches at different times by different authors are gathered. We start with the following example, which considers the following statement *"Every multiplier makes efficient the ST formulas "*

*21°."The generator of line 15 is proposed by G. Marsaglia (G. Marsaglia) as the candidate for the "best multiplier", after computer calculations for almost cubic lattices of dimensionality from 2 to 5. This proposal was made, in particular, because the multiplier can be easily remembered (see the book edited by S.K.Zaremba [Applications of Number Theory to Numerical Analysis, edited by S.K.Zaremba (New York: Academic Press, 1972), 275]. ".*

This generator $\langle X_n(a=69069, N=2^{32})\rangle$ has a maximum period number $N=2^{32}$ as it contains the only prime multiplier, and in $a-1 = d \cdot 2^2 = 17267 \times 2^2$ where $d$ and $N$ are mutually prime, and because the potential is $\tau = \frac{32}{2} = 16$ (see. [2, Section 3.3.1.3, Ex. 5]).

For the orientation in the discussed issue, once we should determine that according to the Theorem 2

$$v_s^2(a,N) \leq a^2 + 1 = 4770526762, \qquad (5.1)$$

for all $s \geq 2$. With this, the data in line 15 of Table 1, beginning with $s=3$, is much less

$$v_3^2 = 2072544, v_4^2 = 52804, v_5^2 = 6990, v_6^2 = 242. \qquad (5.2)$$

Applications of ST-formulas to G.Marsaglia's multiplier $a=69069$ lead to the following results.

For $s=2$ Theorem 1 gives the following exact equality

$$v_2^2(a=69069, N=(69068)^2 = 4770388624) = (a-2)^2 + 1 = (69067)^2 + 1 = 4770250489,$$

that is greater than the value $v_2^2 = 4243209856$ in Table 1 [2, Section 3.3.4] for the value 527040633 calculated for $N = 2^{32} = 4294967296$.

For all $s, s > 2$ the Theorem 7 can be applied,

$$v_s^2(a=69069, N=(69068)^2 = 4770388624) \leq \sum_{k=0}^{2}\left(\binom{2}{k}\right)^2 = 6.$$

Thus for the Marsaglia's multiplier we get the best possible value $v_2^2$ with the maximum period $N = 4770388624$.

To compare the opportunities of ST-formulas with data (5.2) for $\tau = 6$ we apply Theorems 5-6 that for all $s, 2 \leq s \leq 6$ leads to the following lower estimates

$$4767626304 \leq v_s^2(69069; N=(69068)^6). \qquad (5.3)$$

Thus, the random number generator $\langle X_n(69069; N=(69068)^6)\rangle_{n=0}^{(69068)^6-1}$ is constructed with a multiplier of Marsaglia, in which with great maximum period $N = (69068)^6 \approx 2^{192}$ against $N = 2^{32}$ and at the same time large, somewhere with an extremely large (5.1), $v_s^2$ for $s = 3,4,5,6$ in (5.3) against (5.2).

It is obvious that the particular individual properties of Marsaglia's multiplier are not used, since ST-formulas lead to such results for every positive integer $a$.

*22°. «A similar, but less prominent multiplier $16807 = 7^5$ in row 19 became more frequently used for this module after it was offered by Lewis, Goodman and Miller (see Ref. Lewis, Goodman, and Miller in IBM Systems J. 8 (1969), 136-146). The generator with this multiplier is the main factor of a popular library IMSL programs in 1971. The main reason for prolonged use of $a = 16807$ is that $a^2$ is less than module $N$, so the operation $ax \bmod N$ can be performed with high efficiency in high-level languages, using the technique of Ex. 3.2.1.1-9. However, such small multipliers can be contained well known defects. "*[2, Section 3.3.4].

Application of Theorems 1-7, for any $m \geq 2$ (in particular, for $m = 2$ as in considered example) $a - 1 = p_1^{r_1} \cdots p_t^{r_t}$ and $(a-1)^m < (a-1)^{m+t} = N(t \geq 1)$, for example, for $a - 1 = 7^5, N = 7^{10+t}(t \geq 1)$ ensure the fulfillment of $a^m < N$, and with arbitrarily large $a$.